\documentclass[11pt]{article}
\usepackage{amsmath2000,amssymb,amsfonts,amsthm2000,multirow,xspace}
\usepackage[numbers]{natbib}

\makeatletter
\renewcommand{\@seccntformat}[1]{{\csname the#1\endcsname}.\hspace{0.5em}}

\makeatother
\let\section\subsection
\let\subsection\subsubsection

\newtheorem{theorem}{Theorem}
\newtheorem{lemma}[theorem]{Lemma}
\newtheorem{corollary}[theorem]{Corollary}

\newcommand{\sref}[1]{\S~\ref{sec:#1}}
\newcommand{\tref}[1]{Theorem~\ref{thm:#1}}
\newcommand{\lref}[1]{Lemma~\ref{lem:#1}}

\newcommand{\Var}{\operatorname{Var}}
\newcommand{\Cov}{\operatorname{Cov}}

\begin{document}
\title{\vspace*{-62pt}
On the critical exponents of random $k$-SAT}
\author{
\begin{tabular}{c}
David B. Wilson \\
 \small One Microsoft Way\\
 \small Redmond, WA 98052\\
 \small U.S.A.\\
\end{tabular}
}
\date{}
\maketitle

\vspace*{-12pt}
\begin{abstract}
There has been much recent interest in the satisfiability of random
Boolean formulas.  A random $k$-SAT formula is the conjunction of $m$
random clauses, each of which is the disjunction of $k$ literals (a
variable or its negation).  It is known that when the number of
variables $n$ is large, there is a sharp transition from
satisfiability to unsatisfiability; in the case of $2$-SAT this
happens when $m/n\rightarrow 1$, for $3$-SAT the critical ratio is
thought to be $m/n\approx 4.2$.  The sharpness of this transition is
characterized by a critical exponent, sometimes called $\nu=\nu_k$
(the smaller the value of $\nu$ the sharper the transition).
Experiments have suggested that $\nu_3=1.5\pm 0.1$, $\nu_4=1.25\pm0.05$,
$\nu_5=1.1\pm0.05$, $\nu_6=1.05\pm0.05$, and
heuristics have suggested that $\nu_k\rightarrow 1$ as $k\rightarrow\infty$.
We give here a simple proof
that each of these exponents is at least $2$ (provided the exponent is
well-defined).  This result holds for each of the three standard
ensembles of random $k$-SAT formulas: $m$ clauses selected uniformly
at random without replacement, $m$ clauses selected uniformly at
random with replacement, and each clause selected with probability
$p$ independent of the other clauses.  We also obtain similar results
for $q$-colorability and the appearance of a $q$-core in a random
graph.
\end{abstract}

\section{Introduction}

In the past decade many researchers have studied the satisfiability of
random Boolean formulas, in an attempt to understand the ``average
case'' of NP-complete problems.  See \cite{hayes:sat} for a survey.
Let $n$ denote the number of Boolean variables.  A literal is either a
Boolean variable or its negation.  A $k$-clause is the OR
(disjunction) of $k$ literals whose underlying variables are all
distinct.  A random $k$-SAT formula is the AND (conjunction) of $m$
uniformly random $k$-clauses.  A formula is satisfiable if there is an
assignment to the Boolean variables for which the formula evaluates to
TRUE.  For random $3$-SAT it has been observed empirically
\cite{mitchell-selman-levesque} that there
is a critical value $\alpha_3 \approx 4.2$ such that when $n$ is large
and $m/n < \alpha_3 - \varepsilon$, the formula is nearly always
satisfiable, while if $m/n > \alpha_3 + \varepsilon$, the formula is
nearly always unsatisfiable.  Furthermore, determining whether or not
a formula is satisfiable appears to be the hardest when the ratio
$m/n$ is about $\alpha_3$ \cite{mitchell-selman-levesque}.
Similar phenomenona occur for other values $k$, except that for $k=2$
the formulas are always easy (deterministically).  Consequently there
have been many empirical as well as rigorous studies of this
transition from satisfiable to unsatisfiable.

It is known rigorously \cite{friedgut:sharp} that the SAT-to-UNSAT
transition is sharp, i.e.\ that at some critical ratio of $m/n$ the
probability of satisfiability rapidly drops from close to $1$ to close
to $0$.  But for $k>2$ it has not been proved that the critical ratio
of $m/n$ tends to a constant, as opposed to being a slowly varying
function of $n$ that oscillates between its known lower and upper
bounds of $(\log 2) 2^{k-1}-(\log2+1)/2-o_k(1)$ \cite{achlioptas-moore:sat} and
$(\log 2) 2^k$ \cite{chvatal-reed:sat}.  (When $k=3$, the tighter
bounds of $3.42$ \cite{kaporis-kirousis-lalas:3.42} and
$4.506$ \cite{dubois-boufkhad-mandler:3-sat} are known.)

One basic feature of the SAT-to-UNSAT transition is its characteristic
width.  This width is the amount $\Delta$ by which $m$ needs to be
increased for the probability of satisfiability to drop from $2/3$ to
$1/3$, or more generally, to drop from $1-\varepsilon$ to
$\varepsilon$.  The characteristic width is thought to grow as a
polynomial in $n$, so that $\Delta = \Theta(n^{1-1/\nu})$, where the
constant hidden by the $\Theta()$ depends on $\varepsilon$, but the
critical exponent $\nu$ does not.  (Using $\nu$ to denote this
critical exponent is a rather unfortunate choice of notation, since in
statistical mechanics $\nu$ refers to a related but different critical
exponent, see e.g.\ \cite[Chapter 7]{grimmett:percolation}.  The
exponent for the width would be denoted $2-\alpha$, but we use
$\nu$ in this paper to facilitate comparison with earlier studies of
the $k$-SAT transition.)
It is not obvious \textit{a priori} that
$\nu$ is well-defined, as it could in principle slowly oscillate with
$n$ or depend upon $\varepsilon$.  For $2$-SAT it was proved recently
\cite{bollobas-borgs-chayes-kim-wilson} that the characteristic width
does in fact grow polynomially in $n$, and that $\nu=3$.  There have
been a number of experimental studies aimed at measuring the critical
exponent $\nu$ for random $k$-SAT, as summarized in the above table,
and several authors have conjectured that the exponent $\nu$ tends to
1 as $k$ gets large.  The purpose of this note is to provide a simple
proof that for each fixed $k$, the characteristic width is always at
least $\Theta(n^{1/2})$, so that in particular, if the exponent $\nu$
is well-defined, it is always at least $2$.

\newcommand{\citen}[2]{\citeauthor*{#2} \citep*[#1]{#2}}
\renewcommand{\thefootnote}{\fnsymbol{footnote}}
\newcommand{\tfn}[2]{\protect\newline{\footnotesize\footnotemark[#1] #2}}
\newcommand{\tfnm}[1]{\protect\footnotemark[#1]}
\begin{table}
\renewcommand{\multirowsetup}{\centering}
\newlength{\KK}\settowidth{\KK}{$k$}
\centerline{
\parbox{5in}{
\begin{tabular}{|c||c|c|}%|c|c|}
\hline
\multirow{2}{\KK}{$k$} & \multicolumn{2}{|c|}{``$\nu$'' for SAT-to-UNSAT transition width}  \\%& \multicolumn{2}{|c|}{``$\nu$'' for computational cost} \\
\cline{2-3}
    & conjectured & rigorous  \\% & conjectured & rigorous \\
\hline
\hline
$2$ & $2.6\pm 0.2$ \tfnm{1} \tfnm{2}, $2.8$ \tfnm{3}, $3$ \tfnm{8} \tfnm{6} & $3$ \tfnm{9} \\% &&\\
\hline
$3$ & $1.5\pm 0.1$ \tfnm{1} \tfnm{2}, $1.5$ \tfnm{3} \tfnm{4} & $\geq 2$ \tfnm{7} \\% & $1.3$ \tfnm{8} & $\geq 2$ \tfnm{7} \\
\hline
$4$ & $1.25\pm 0.05$ \tfnm{1} \tfnm{2}, $1.25$ \tfnm{5} & $\geq 2$ \tfnm{7} \\% & $1.25$ \tfnm{8} & $\geq 2$ \tfnm{7} \\
\hline
$5$ & $1.10\pm 0.05$ \tfnm{1} \tfnm{2} & $\geq 2$ \tfnm{7} \\% & $1.1$ \tfnm{8} & $\geq 2$ \tfnm{7} \\
\hline
$6$ & $1.05\pm 0.05$ \tfnm{1} \tfnm{2} & $\geq 2$ \tfnm{7} \\% & & $\geq 2$ \tfnm{7} \\
\hline
$\rightarrow\infty$ & $\rightarrow 1$ \tfnm{1} \tfnm{2} \tfnm{3} & $\geq 2$ \tfnm{7} \\% & & $\geq 2$ \tfnm{7} \\
\hline
\end{tabular}
\smallskip
\tfn{1}{\citet*{kirkpatrick-selman:critical}}
\tfn{2}{\citet*{kirkpatrick-gyorgyi-tishby-troyansky:k-sat}}
\tfn{3}{\citet*{monasson-zecchina-kirkpatrick-selman-troyansky:rsa} page 428}
\tfn{4}{\citet*{monasson-zecchina-kirkpatrick-selman-troyansky:nature} page 135}
\tfn{6}{\citet*{monasson-zecchina:replica}}
\tfn{5}{\citet*{gent-walsh:sat} page 109}
\tfn{8}{\citet*{bollobas-borgs-chayes-kim}}
\tfn{9}{\citet*{bollobas-borgs-chayes-kim-wilson}}
\tfn{7}{This article.}
}
}%\end{center}
\end{table}

\textit{Remark:\/}
There is a related ensemble of random $k$-SAT formulas, in which each
possible $k$-clause appears in the formula with probability $p$
independently of the other clauses.  For convenience we let $M$ denote
the total number $2^k\binom{n}{k}$ of possible clauses.  When $p M
\approx m$, this ${\cal F}_{n,p}$ ensemble of random formulas will
behave much like the ${\cal F}_{n,m}$ ensemble of formulas defined
above.  But there is a limit on how sharp the SAT-to-UNSAT transition
can be for the ${\cal F}_{n,p}$ ensemble, due to the approximate
relationship between $p$ and the number of clauses in the formula.
Even if $p M = m$, the number of clauses will be $m \pm
\Theta(m^{1/2})$.  It is thus straightforward to show
that in the ${\cal F}_{n,p}$ ensemble, the critical exponent $\nu$
must be at least $2$.  (This bound is closely related to the ``Harris
criterion'' for disordered statistical mechanical systems
\cite{chayes-chayes-fisher-spencer:prl}
\cite{chayes-chayes-fisher-spencer:cmp}.)  The (by now) standard proof of
$\nu\geq 2$ for the ${\cal F}_{n,p}$ ensemble depends on the variance
in the number of clauses, which is zero for the ${\cal F}_{n,m}$
ensemble.  It is simple to define properties on sets of clauses
such that there is a much sharper transition in the ${\cal F}_{n,m}$
ensemble, with a smaller value of $\nu<2$ --- one trivial example is
the property ``$m\leq 4.2 n$'', for which $\nu=1$ in ${\cal F}_{n,m}$
while $\nu=2$ in ${\cal F}_{n,p}$.
Until now it has been
suggested (on the strength of Monte Carlo experiments and
heuristic arguments) that satisfiability is one such property.
Despite this, we proceed to show that even for the ${\cal F}_{n,m}$
ensemble, the characteristic width is at least $\Theta(n^{1/2})$, so
that $\nu$ (if it is well-defined) must be at least $2$.

There are also a number of questionable
conjectures about other features of random $k$-SAT formulas.  For
instance, for $3$-SAT \citet*{crawford-auton:3-sat} study the value of
$m=m_{1/2}(n)$, where $m_r(n)$ denotes the smallest value of $m$ for
which the fraction of satisfiable $k$-SAT formulas is $\leq r$.  They
fit $m_{1/2}(n)$ to a curve of the form $m_{1/2}(n)=\alpha_3 n + A
n^{1-1/\nu}$, obtaining $\nu=1$ \citep*{crawford-auton:sat} and then
later $\nu=3/5$ \citep*{crawford-auton:3-sat}.  (The integrality of $m$
by itself strongly suggests that $\nu\geq 1$.)
While it is in principle conceivable that $m_{1/2}(n)=\alpha_3 n + O(1)$,
it is an easy consequence of our work that there can be at most one
value of $r$ for which $m_r(n) = \alpha_3 n + o(n^{1/2})$ (otherwise
the SAT-to-UNSAT transition would be too sharp).  For $2$-SAT the
special value of $r$ is empirically about $91\%$, and for $3$-SAT
there is no reason to believe (and experimental reason not to believe)
that the special value of $r$ is $1/2$.

\citet*{selman-kirkpatrick:cost} estimated exponents for the
characteristic width of the median computational difficulty of
determining whether or not a random formula is satisfiable, where
computational difficulty was measured in terms of the number of
recursive calls made by Crawford and Auton's SAT-solver
(\textsc{Tableau}).  The exponents they obtained were $1.3$ for
$3$-SAT, $1.25$ for $4$-SAT, and $1.1$ for $5$-SAT.  We do not analyze
the specific SAT-solver \textsc{Tableau}, but we can say something
about the characteristic width of the computationally difficult
problems for other SAT-solvers.  Many SAT-solvers use the ``pure
literal rule'' before starting a backtracking search for a satisfying
assignment.  A SAT-solver using this rule will look for literals $y$
in the formula such that $y$'s negation $\bar{y}$ does not also appear
in the formula.  If such a literal $y$ exists, then the SAT-solver
sets $y$ to TRUE and removes from the formula any clauses containing
$y$, since the resulting simpler formula is satisfiable if and only if
the original formula was satisfiable.  Rigorously analyzing the
\textit{median\/} computational difficulty seems not so easy, but
using our methods one can show that for these SAT-solvers the
\textit{typical\/} computational difficulty has critical exponent at least $2$.
By this we mean the following: if for $m$ clauses there is probability
$p$ that the number of recursive calls is between $L$ and $U$, then
when there are $m+\Delta$ clauses, the probability is
$p+O(\Delta/\sqrt{n})+o(1)$.

Our method is general enough to be applicable to other types of sharp
transitions.  For instance, \citet*{pittel-spencer-wormald:k-core}
prove that there is a sharp transition for the appearance of a
$q$-core in a random graph.  (The $q$-core is the maximal subgraph for
which each vertex has degree at least $q$.)  They prove that the width
of this transition is at most $n^{1/2+o(1)}$, but gave no lower bound.
We prove that the width is at least $\Theta(n^{1/2})$.  (Independent
of this present work, \citet*{kirkpatrick:personal} has reported that
experiments suggest that the width is $\Theta(n^{1/2})$.)
We can also supply
lower bounds on the transition width of other graph properties, such
as $q$-colorability.

\section{Proofs of theorems}
\label{sec:theorems}

We now give a proof that $\nu\geq 2$ in ${\cal F}_{n,m}$ that works
simultaneously for $k$-SAT, $q$-colorability, the existence of a
$q$-core, and a variety of other properties.
Let $M$ denote the number of possible items.  In the case of $k$-SAT, the items
will be the possible clauses on $n$ variables, and
$M=2^k\binom{n}{k}$.  In the case of the $q$-core or $q$-colorability,
the items will be the possible edges of a graph on $n$ vertices, and
$M=\binom{n}{2}$.  A \textit{property\/} classifies sets of items into
two types: sets which have the property (for convenience call them
\textit{proper\/} sets), and sets which do not have the property
(\textit{improper\/} sets).
Satisfiability is a property on sets of clauses, $q$-colorability is a
property on sets of edges, and the existence of a $q$-core of a graph
is also a property on sets of edges.  We may also be interested in
non-monotone properties, such as the property that a certain
SAT-solver does more than $L$ recursive calls to determine the
satisfiability of a set of clauses.

A \textit{bystander rule\/} is a way of partitioning a set of items into
two classes: the relevant items, and the bystander items.  A bystander
rule must satisfy the following constraint.  Given a set of items $A$,
let $R$ be the relevant items, and let $G$ be the bystander items.  Then
for any $B\subseteq A$, it must be that $B$ has the property if and
only if $B\setminus G$ has the property.  In this sense, the only relevant
items for the property are those that are contained in $R$ --- that is, if
one restricts attention to sets of items contained in $A$,
the bystander items never affect the property.

In the case of $k$-SAT, we will use the ``partially-free'' bystander rule,
which declares a clause to be a bystander if the underlying variable of
one of its literals does not appear anywhere else in the formula.
(Recall that the formula is the AND of the clauses in the set.)  By
setting this variable to an appropriate value, the clause can be
satisfied without affecting our ability to satisfy the remaining
clauses of the formula.  Thus partially-free is in fact a valid bystander
rule.

For the $q$-core and $q$-colorability, we also use the
``partially-free'' bystander rule, which in the context of graphs (and
hypergraphs) declares an edge to be a bystander if one of its endpoints
has degree $1$.  It is simple to check that partially-free is a valid
bystander rule for these properties.

\textit{Remark:\/}
We could in principle use instead other bystander rules.  One
possibility is to declare any clause that eventually gets
resolved by repeated application of the pure literal
rule to be a bystander.  But it is easier to analyze the partially-free bystander rule,
and our objective to provide a \textit{simple\/} rigorous proof.

\begin{theorem}
\label{thm:sqrt}
Suppose that a property has a bystander rule such that, in a set of $m$
random items, with probability $1-\varepsilon$ at least $\gamma m$ of
the items are bystanders.  (Items may be chosen either without or with
replacement.)  Suppose further that a set of $m_1<m$ random
items is proper with probability $p_1$, and a set of $m_2<m$ random
items is proper with probability $p_2$.
If $\beta\leq m_1/m\leq 1-\beta$, and $\beta\leq m_2/m\leq 1-\beta$, then
$$|m_1 - m_2| \geq \left(|p_1-p_2|-\varepsilon\right) \sqrt{2\pi m}
 \sqrt{\frac{\gamma \beta (1-\beta)}{1-\gamma}}(1-o(1)),$$
where the $o(1)$ term becomes small if $m$ gets large while $\beta$
and $\gamma$ remain fixed.
\end{theorem}

Informally, \tref{sqrt} says that if there are many bystanders, then
the transition cannot be too sharp.  To prove this we use the following lemma:
\begin{lemma}
\label{lem:rg-balls}
Suppose there are $m$ balls, of which $\gamma m$ are green and $(1-\gamma) m$ are red, and that $\beta m$ of these
balls are randomly sampled without replacement.  The probability that
the $(\beta m)$th ball is red and exactly $\ell$ of the sampled balls are red
is, as a function of $\ell$, is unimodal and at most
$$ \frac{1+o(1)}{\sqrt{2\pi m}}\sqrt{\frac{1-\gamma}{\gamma\beta(1-\beta)}}.$$
Here the $o(1)$ term becomes small when the expected number of sampled
red balls, sampled green balls, unsampled red balls, and unsampled
green balls are each large.
\end{lemma}
\noindent
We postpone the proof of this lemma to \sref{lemmas}, and proceed to the more interesting
part of the proof.

\begin{proof}[Proof of \tref{sqrt}]
Let $C_1,\ldots,C_M$ denote the items.  If the items are selected
without replacement, let $\sigma$ be a uniformly random permutation on
the numbers $1,\ldots,M$.  If the items are selected with
replacement, let $\sigma$ be an i.i.d.\ sequence of uniformly random
integers in the range $1,\ldots,M$.  Let $f_m$ denote
the sequence consisting of the first $m$ items with respect to $\sigma$,
i.e.\ $f_m = \langle C_{\sigma(1)},\ldots,C_{\sigma(m)}\rangle$; $f_m$ is a
uniformly random sequence of $m$ items chosen without (resp.\ with) replacement.
Let $g_m$ be the number of bystander
items, and $r_m=m-g_m$ the number of relevant items of $f_m$.  Say that $\ell$
is a positively (resp.\ negatively) critical integer if the set of the
first $\ell$ relevant items of $f_m$ is proper
(resp.\ improper) but the first $\ell-1$ relevant items is improper
(resp.\ proper).  Pick a uniformly random permutation $\tau$ on the
numbers $1,\ldots,m$ (independent of $\sigma$).  Use $\tau$ to tag a
random set of $b$ items from $f_m$, i.e.\
$C_{\sigma(\tau(1))},\ldots,C_{\sigma(\tau(b))}$; the tagged items
form a random set of $b$ items chosen without (resp.\ with) replacement,
since we could have picked $\tau$
first and then $\sigma$.  Suppose that $L$ of the tagged items are
relevant.  Since our sequence of items $f_m$ is already in a random order,
we may instead pick and keep the first $L$ relevant items of $f_m$
and the first $b-L$ bystander items of $f_m$.  The resulting set
$\hat{f}_{m,b}$ of kept items is a uniformly random set of $b$
items chosen without (resp.\ with) replacement, and whether
or not $\hat{f}_{m,b}$ has the property is determined by $L$.  We can write
\begin{align*}
\Pr[\text{$\hat{f}_{m,b}$ is proper}|f_m] - \Pr[\text{$\hat{f}_{m,b-1}$ is proper}|f_m] =&
  \Pr[\text{$\hat{f}_{m,b}$ is proper and $\hat{f}_{m,b-1}$ is improper}|f_m]\\
&-\Pr[\text{$\hat{f}_{m,b}$ is improper and $\hat{f}_{m,b-1}$ is proper}|f_m]\\
=&\Pr\left[\parbox{2.47in}{\# relevant tags is positively critical, $b$th tag is relevant}\big|f_m\right]\\
&-\Pr\left[\parbox{2.47in}{\# relevant tags is negatively critical, $b$th tag is relevant}\big|f_m\right]
\end{align*}
Now we use the fact that the negatively critical integers and the
positively critical integers are interleaved, and that
$$f(\ell)=\Pr[\text{\# relevant tags is $\ell$, $b$th tag is relevant}|f_m]$$
is unimodal in $\ell$ (from \lref{rg-balls}).  Let the critical integers be
$\ell_1,\ell_2,\ldots,\ell_c$, and suppose that of these, $\ell_\mu$
maximizes $f()$.  Then we can write
\begin{align*}
 \sum_{i=1}^c (-1)^{i-\mu} f(\ell_i) &= f(\ell_\mu) \begin{aligned}[t]
 -& [f(\ell_{\mu+1})-f(\ell_{\mu+2})] - \cdots\\
 -& [f(\ell_{\mu-1})-f(\ell_{\mu-2})] - \cdots
\end{aligned}\\
&\leq f(\ell_\mu)
\end{align*}
and
\begin{align*}
 \sum_{i=1}^c (-1)^{i-\mu} f(\ell_i) &= f(\ell_\mu) \begin{aligned}[t]
 -& f(\ell_{\mu+1}) + [f(\ell_{\mu+2})-f(\ell_{\mu+3})] + \cdots\\
 -& f(\ell_{\mu-1}) + [f(\ell_{\mu-2})-f(\ell_{\mu-3})] + \cdots
 \end{aligned} \\
&\geq f(\ell_\mu) - f(\ell_{\mu+1}) - f(\ell_{\mu-1}) \geq -f(\ell_\mu).
\end{align*}
Thus
\begin{align*}
\Big|\Pr[\text{$\hat{f}_{m,b}$ is proper}|f_m] - \Pr[\text{$\hat{f}_{m,b-1}$ is proper}|f_m]\Big| 
&\leq \max_\ell \Pr\left[\parbox{1.3in}{\# relevant tags is $\ell$,  $b$th tag is relevant}\Big|f_m\right] \\
&\leq \frac{1+o(1)}{\sqrt{2\pi m}} \sqrt{\frac{r_m/m}{(g_m/m) (b/m) (1-b/m)}}
\intertext{by \lref{rg-balls}, and then assuming $g_m \geq \gamma m$ and $\beta m \leq b \leq (1-\beta)m$ we get}
&\leq \frac{1+o(1)}{\sqrt{2\pi m}} \sqrt{\frac{1-\gamma}{\gamma \beta (1-\beta)}},
\end{align*}
where the $o(1)$ term becomes small if $m$
gets large while $\beta$ and $\gamma$ remain fixed.  Thus if both $m_1$ and $m_2$ are between $\beta m$ and $(1-\beta m)$ we can write
\begin{align*}
\left|
\begin{aligned}\ &\Pr[\text{$\hat{f}_{m,m_1}$ is proper}]\\
              -&\Pr[\text{$\hat{f}_{m,m_2}$ is proper}]
\end{aligned}
\right|
  \leq&
\begin{aligned}
\ &\Pr[\text{$f_m$ has at least $\gamma m$ bystanders}] \times |m_1-m_2|
  \frac{1+o(1)}{\sqrt{2\pi m}} \sqrt{\frac{1-\gamma}{\gamma \beta (1-\beta)}}\\
+ &\Pr[\text{$f_m$ has less than $\gamma m$ bystanders}]
\end{aligned}\\
|p_1-p_2|
\leq& |m_1-m_2| \frac{1+o(1)}{\sqrt{2\pi m}} \sqrt{\frac{1-\gamma}{\gamma \beta (1-\beta)}} + \varepsilon. \qedhere
\end{align*}
\end{proof}

To apply \tref{sqrt} to $k$-SAT, we need
to show that many clauses are bystanders, and to apply it to the $q$-core
and $q$-colorability thresholds, we need to show that many edges are
bystanders.

\begin{lemma}
\label{lem:many-free}
Suppose that the items are $k$-clauses, edges, or hyperedges on $k$
vertices.  Assume that $k$ is fixed, and that $m=O(n)$ random items
are selected, either with replacement or without replacement.  With
high probability there will be $(1+o(1)) m [1-[1-e^{-km/n}]^k]$
partially-free items.
\end{lemma}

This lemma says that the number of partially-free clauses or hyperedges
is with high probability close to what one might naively expect.
Since we use the lemma to disprove a number of experimental results
and heuristic arguments, we give a careful proof of it in \sref{lemmas}.
But first let us see how to use the lemma with \tref{sqrt}.

\begin{corollary}
\label{thm:3-sat}
Let $p_1$ and $p_2$ be fixed numbers such that $1>p_1>p_2>0$.
Suppose that a random $3$-SAT formula with $n$ variables and $m_1$ clauses is
satisfiable with probability $\geq p_1$ and a random $3$-SAT formula
with $n$ variables and $m_2$ clauses is satisfiable with probability $\leq p_2$.
Then $$m_2-m_1 \geq (0.0015 +o(1)) \times (p_1-p_2) \times \sqrt{n},$$
where the $o(1)$ goes to $0$ when $n\rightarrow\infty$.  In
particular, $\nu_3\geq 2$ if it is well-defined.  Similarly, for
$k$-SAT in general ($k$ fixed), we get $m_2-m_1 \geq
(p_1-p_2)\Theta(\sqrt{n})$, implying $\nu_k\geq 2$.
\end{corollary}

\begin{proof}
For $3$-SAT, when $n$ is large, we know that $m_1/n$ and $m_2/n$ are
both close to the critical ratio $c_3(n)$, where $3.42\leq c_3(n)
\leq 4.571$.  (Details of the $4.506$ upper bound are not available
at this time, so here we use the established bound of $4.571$
\cite{kaporis-kirousis-stamatiou-vamvakari-zito:sat}.)
More generally for $k$-SAT, we know that for $i=1$ or $2$,
$\check{c}_k-o(1)\leq m_i/n \leq \hat{c}_k+o(1)$, where $\check{c}_k$ and
$\hat{c}_k$ are lower and upper bounds on the critical $k$-SAT ratio
$c_k(n)$ (which could conceivably be a function of $n$).  Let
$m=n(c_k(n) + t)$, where $t$ is a positive constant that
we will choose in a moment.  By \lref{many-free}, the fraction
of clauses which are partially-free is w.h.p.\
$\gamma = 1-[1-e^{-k(c_k(n)+t)}]^k+o(1)$.  The value $\beta$ for
\tref{sqrt} is $t/(c_k(n)+t)$.  Note that
$\gamma$ is monotone decreasing in $c_k(n)$, and that
$m \beta(1-\beta) = n t c_k(n) / (c_k(n)+t)$ is monotone increasing in
$c_k(n)$, so that our bound from \tref{sqrt} is at least as good as
\begin{align*}
|m_1 - m_2|
&\geq \left(|p_1-p_2|-o(1)\right) \sqrt{2\pi n}
\sqrt{\frac{\check{c}_k t}{\check{c}_k+t}}
\sqrt{\frac{1-[1-e^{-k(\hat{c}_k+t)}]^k}{[1-e^{-k(\hat{c}_k+t)}]^k}}(1-o(1)) \\
&= (p_1-p_2)\Theta(\sqrt{n}).
\end{align*}
For $3$-SAT we take $t=0.3$ to get the above-stated constant of $0.0015$.
\end{proof}

For the appearance of a $q$-core in a random graph, there is a sharp
threshold, and furthermore the precise values of the critical ratio
$m/n = c_q$ are known.  It is known e.g.\ that $c_3\approx 3.35$,
$c_4\approx 5.14$, $c_5\approx 6.81$, and $c_k=k+\sqrt{k\log k}+O(\log
k)$ \cite{pittel-spencer-wormald:k-core}.  We can lower bound the
characteristic width for the appearance of the $q$-core in essentially
the same that we did for $k$-SAT, except that here $k=2$ even as $q$
varies.  (A larger value of $k$ would correspond to the appearance of
the $q$-core within a $k$-uniform hypergraph.)

\begin{corollary}
\label{thm:q-core}
For $q\geq3$, the transition for existence of a $q$-core has characteristic width $\geq\Theta(\sqrt{n})$.
\end{corollary}

\noindent
When one randomly adds edges one a time, w.h.p.\ the $q$-core jumps
from size $0$ to size $\Theta(n)$ with the addition of a single edge
\cite{pittel-spencer-wormald:k-core} ---
Corollary~\ref{thm:q-core} is a statement about the timing of this jump.

With $q$-colorability, it is known that there is a sharp threshold
\cite{achlioptas-friedgut:sharp} for the number of edges that a random
graph can have while still being $q$-colorable, but as with $k$-SAT,
it is not known that $c_q$ is a \textit{bona fide\/} constant rather
than a slowly varying function of $n$ that oscillates between its
known upper and lower bounds.  \L uczak proved
\nocite{luczak:chromatic} proved that $c_q/(q\log q) \rightarrow 1$ as
$q\rightarrow\infty$, and it is known that $2.01 \leq c_3 \leq 2.495$
\cite{achlioptas-moore:coloring} \cite{kaporis-kirousis-stamatiou:coloring}.
As with the $q$-core, here $k=2$ even as $q$ varies, and we have

\begin{corollary}
\label{thm:q-color}
For $q\geq 3$, the transition for $q$-colorability has characteristic
width at least $\Theta(\sqrt{n})$.
\end{corollary}

\section{Proofs of lemmas}
\label{sec:lemmas}

\begin{proof}[Proof of \lref{rg-balls}]
For convenience let $g=\gamma m$ be the number of green balls,
$r=(1-\gamma)m$ be the number of red balls, and $b=\beta m$ be the
number of balls selected.
The precise probability that the $b$th ball is the $\ell$th red ball is
$$ \frac{\binom{r}{\ell} \binom{g}{b-\ell} \frac{\ell}{b}}{\binom{r+g}{b}} .$$
The ratio of successive terms is
\begin{equation*}
\frac{\binom{r}{\ell} \binom{g}{b-\ell} \frac{\ell}{b}}{\binom{r+g}{b}}
\frac{\binom{r+g}{b}}{\binom{r}{\ell+1} \binom{g}{b-\ell-1} \frac{\ell+1}{b}}
 =
\frac{\ell(g-b+\ell+1)}{(r-\ell) (b-\ell)}
\end{equation*}
which is monotone in $\ell$ and implies the unimodality claim.  Next we
identify the mode:
\begin{align*}
\ell(g-b+\ell+1) &\geq (r-\ell) (b-\ell) \\
\ell(g-b+1+b+r) &\geq r b \\
\ell &\geq \frac{r b}{r+g+1}
\end{align*}
so that the optimal $\ell$ is given by
$\ell = \lceil r b/(r+g+1)\rceil$, which is within $1$ of $r b /(r+g)$.

We next approximate the maximum value of this probability.  For
convenience let $\lambda = \ell/m$.
Recall Stirling's formula: $n! = n^n e^{-n} \sqrt{2\pi n}
\exp[1/(12n+\delta_n)]$ where $0\leq \delta_n \leq 1$.
\begin{align*}
\frac{\binom{r}{\ell} \binom{g}{b-\ell} \frac{\ell}{b}}{\binom{m}{b}}
=& \frac{\ell}{b} \times
\frac{((1-\gamma) m)!(\gamma m)!(\beta m)!((1-\beta)m)!}
     {(\lambda m)!((1-\gamma-\lambda)m)!
      ((\beta-\lambda)m)!((\gamma-\beta+\lambda)m)! m!} \\
=& \frac{\lambda}{\beta} \times
\left[
\frac{(1-\gamma)^{1-\gamma}\gamma^{\gamma}\beta^{\beta}(1-\beta)^{1-\beta}}
     {(\lambda)^{\lambda}(1-\gamma-\lambda)^{1-\gamma-\lambda}
      (\beta-\lambda)^{\beta-\lambda}
      (\gamma-\beta+\lambda)^{\gamma-\beta+\lambda}} 
\right]^m \\
& \times
\frac{1}{\sqrt{2\pi m}}
\sqrt{
\frac{(1-\gamma) \gamma\beta (1-\beta)}
     {\lambda (1-\gamma-\lambda)
      (\beta-\lambda)(\gamma-\beta+\lambda)}}
\times \exp(o(1)),
\end{align*}
Consider the $\exp(o(1))$ error term arising from the $\exp[1/(12
n+\delta_n)]$ portion of Stirling's formula.  Since $\ell\leq r$,
$r-\ell\leq m-b$, $b-\ell\leq b$, and $g-b+\ell\leq g$, the error term
will be $\leq 1$, so we may drop it to get an upper bound.
If the second term on the right were larger than $1$, then we could
increase $m$ while keeping the ratios $g/m$, $\ell/m$, and $b/m$
fixed, and thereby make the probability as large as we like, and in
particular larger than $1$.  Thus we can drop this term as well:
\begin{align*}
\frac{\binom{r}{\ell} \binom{g}{b-\ell} \frac{\ell}{b}}{\binom{m}{b}}
 &\leq 
\frac{1}{\sqrt{2\pi m}}
\sqrt{
\frac{(1-\gamma) \gamma\lambda (1-\beta)}
     {\beta (1-\gamma-\lambda)
      (\beta-\lambda)(\gamma-\beta+\lambda)}}\\
\intertext{and upon substituting $\lambda=(1-\gamma)\beta \pm 1/m$ we find}
 &\leq 
\frac{1}{\sqrt{2\pi m}}
\sqrt{
\frac{(1-\gamma) \gamma(1-\gamma)\beta (1-\beta)}
     {\beta ((1-\gamma)(1-\beta))
      (\gamma\beta)(\gamma(1-\beta))}}\times (1+o(1))\\
 &\leq 
\frac{1}{\sqrt{2\pi m}}
\sqrt{
\frac{1-\gamma }
     {\beta(1-\beta)\gamma}}\times (1+o(1)),
\end{align*}
where the $o(1)$ vanishes when $\beta(1-\gamma)\gg 1/m$, $\beta\gamma\gg
1/m$, $(1-\beta)(1-\gamma)\gg 1/m$, and $(1-\beta)\gamma\gg 1/m$.
\end{proof}

\textit{Remark:\/}
One referee suggested an alternative to \lref{many-free}, which has
worse constants but a much shorter proof.  The proof of the
alternative lemma focuses on ``nearly free'' variables/vertices rather
than clauses/edges, and uses standard large deviation inequalities for
martingales.  We give here the original proof since it requires less
background.

\begin{proof}[Proof of \lref{many-free}]
Let $r=2$ if the items we are interested in are clauses of $k$ Boolean
variables, and let $r=1$ if the items are edges of a graph ($k=2$) or
$k$-uniform hypergraph ($k\geq 2$).  Recall that $n$ is the number of
variables or vertices.  In each of these cases, the number $M$ of
possible items is $M=r^k\binom{n}{k}$.  Recall our assumption that $k$
is fixed and that we are looking at sets of $m=O(n)$ items, so that
functions of $k$ and $m/n$ may be written as $O(1)$.  A more
detailed analysis could determine what happens when e.g.\
$k\rightarrow\infty$ with $n$, but we do not attempt this.

Say that an item is $d$-free ($1\leq d\leq k$), with respect to a set
of items, if the first $d$ variables (if the item is a clause) or the
first $d$ vertices (if the item is an edge or hyperedge) do not occur
in any other items in the set.

The probability that an item is $d$-free when the $m$ items are randomly
selected without replacement is easily seen to be
\begin{equation*}
\Pr[\text{item is $d$-free}] 
= \frac{\binom{r^k\binom{n-d}{k}}{m-1}}
        {\binom{r^k\binom{n}{k}-1}{m-1}}
= \prod_{i=0}^{m-2} 
\frac{r^k \frac{(n-d)(n-d-1)\cdots (n-d-k+1)}{k!} - i}
     {r^k \frac{(n)(n-1)\cdots (n-k+1)}{k!} - i - 1}.
\end{equation*}
We use the identity
$$ \frac{a-\delta}{b-\delta} = \frac{a}{b} - \frac{1-a/b}{b-\delta}\delta$$
to estimate each term in the product:
\begin{align*}
\frac{r^k \frac{(n-d)_k}{k!} - i}
     {r^k \frac{(n)_k}{k!} - i - 1}
&= 
\frac{(n-d)_k}{(n)_k} - \frac{1-(n-d)_k/(n)_k}{M-(i+1)}(i+1) + \frac{1}{M-(i+1)} \\
&=
\frac{(n-d)_k}{(n)_k} - \frac{(1+o(1))dk/n}{M}O(m) + O(1/M) \\
&=
\frac{(n-d)_k}{(n)_k} \exp[O(1/M)] = \exp[-dk/n+O(1/n^2)]
\end{align*}
(When sampling is done with replacement, the $\exp[O(1/M)]$ error term
does not appear.)
From this we see that the probability that an item in the set is $d$-free is $\exp[-(1+o(1)) d k m/n]$.
Thus the expected number of $d$-free items is $m \exp[-(1+o(1)) dkm/n]$.  We wish to show that the actual number of $d$-free items will likely be close to its expected value, so we bound the variance.

Let $X_C^{(d)}$ be the indicator random variable for item $C$ being $d$-free.  For $C'\neq C$, when sampling is done without replacement we have
\begin{align*}
\Pr[\text{items $C$ and $C'$ $d$-free}]
&= \Pr\left[\parbox{2in}{first $d$ variables/vertices of $C$ are not in $C'$ \& vice versa}\right] \times
   \frac{\binom{r^k\binom{n-2d}{k}}{m-2}}
        {\binom{r^k\binom{n}{k}-2}{m-2}} \\
&= \left(1-\frac{2dk-d^2+o(1)}{n}\right) \times
 \prod_{i=0}^{m-3} 
 \frac{r^k \frac{(n-2d)(n-2d-1)\cdots (n-2d-k+1)}{k!} - i}
     {r^k \frac{(n)(n-1)\cdots (n-k+1)}{k!} - i - 2}.
\end{align*}
In the same manner as above we estimate each term in the product:
\begin{align*}
 \frac{r^k \frac{(n-2d)_k}{k!} - i}
     {r^k \frac{(n)_k}{k!} - i - 2}
&=
\frac{(n-2d)_k}{(n)_k} \exp[O(1/M)].
\end{align*}
(As before, when sampling with replacement, the $\exp[O(1/M)]$
error term does not appear.)
Thus
\begin{align*}
\frac{E[X_C^{(d)} X_{C'}^{(d)}]}{E[X_C^{(d)}] E[X_{C'}^{(d)}]}
&=\frac{1-(2kd-d^2+o(1))/n}{(1-(kd+o(1))/n)^2}\left(\frac{(n-2d)_k}{(n)_k}\frac{(n)_k}{(n-d)_k}\frac{(n)_k}{(n-d)_k}\right)^{m-2} \exp(O(m/n^k))\\
&=\left(1+\frac{d^2+o(1)}{n}\right)\left(\prod_{j=0}^{k-1} \frac{(n-2d-j)(n-j)}{(n-d-j)(n-d-j)}\right)^{m-2}\exp(O(m/n^k))\\
&=\left(1+\frac{d^2+o(1)}{n}\right)\left(\prod_{j=0}^{k-1} \left(1-(1+o(1)) \frac{d^2}{n^2}\right)\right)^{m-2}\exp(O(m/n^k))\\
&= 1 + O(1/n) \\
E[X_C^{(d)} X_{C'}^{(d)}] &= E[X_C^{(d)}] E[X_{C'}^{(d)}] + O(1/n) \\
\Cov\left(X_C^{(d)},X_{C'}^{(d)}\right) &= O(1/n),
\end{align*}
yielding the variance in the number of $d$-free items
\begin{align*}
\Var\left[\sum_C X_C^{(d)}\right] =& m(m-1) \Cov\left(X_C^{(d)},X_{C'}^{(d)}\right) + m \Var\left(X_C^{(d)}\right) = O(m).
\end{align*}
Using Chebychev's inequality, it follows that the actual number of $d$-free items will with high probability be within $O(\sqrt{m})$ of its expected value $m \exp[-(1+o(1))dkm/n]$.

Recall that an item is partially free if at least one of its
variables/vertices is not contained in any of the other items.
We use inclusion-exclusion to estimate the number of partially
free items.  Let $X_C^+$ denote the event that item $C$ is partially
free, and $X_C^{j_1,\ldots,j_d}$ denote the event that item $C$ is
free in positions $j_1,\ldots,j_d$.  Then
\begin{align*}
X_C^+ &= 1-\sum_{S\subseteq \{1,\ldots,k\}} (-1)^{\# S} X_C^S \\
\sum_C X_C^+ &= m-\sum_{S\subseteq \{1,\ldots,k\}} (-1)^{\# S} \sum_C X_C^S.
\end{align*}
Since there are $2^k = O(1)$ possible values of $S$, and for each one
$\sum_C X_C^S$ is with high probability within $O(\sqrt{m})$ of 
$m \exp[-(1+o(1))(\#S)km/n]$, the number of partially free clauses is
with high probability within $O(\sqrt{m})$ of
\begin{align*}
 m-\sum_{S\subseteq \{1,\ldots,k\}} (-1)^{\# S} m \exp[-(1+o(1))(\#S)km/n]
=& m(1+o(1)) \left[1-\sum_{d=0}^k \binom{n}{k} (-1)^d e^{-dkm/n}\right] \\
=& m(1+o(1)) [1-[1-e^{-km/n}]^k]. \qedhere
\end{align*}
\end{proof}

\section*{Acknowledgements:} We thank the referees for their comments.

\newpage

\bibliography{sqrt}

\begin{thebibliography}{28}
\expandafter\ifx\csname natexlab\endcsname\relax\def\natexlab#1{#1}\fi

\bibitem[Achlioptas and Friedgut(1999)]{achlioptas-friedgut:sharp}
Dimitris Achlioptas and Ehud Friedgut.
\newblock A sharp threshold for $k$-colorability.
\newblock {\em Random Structures \& Algorithms}, 14\penalty0 (1):\penalty0
  63--70, 1999.

\bibitem[Achlioptas and Moore(2002{\natexlab{a}})]{achlioptas-moore:coloring}
Dimitris Achlioptas and Cristopher Moore.
\newblock Almost all graphs with average degree 4 are 3-colorable.
\newblock In {\em Proceedings of the Thirty-Fourth Annual ACM Symposium on
  Theory of Computing}, pages 199--208, 2002{\natexlab{a}}.

\bibitem[Achlioptas and Moore(2002{\natexlab{b}})]{achlioptas-moore:sat}
Dimitris Achlioptas and Cristopher Moore.
\newblock The asymptotic order of the $k$-{SAT} threshold, 2002{\natexlab{b}}.
\newblock Manuscript.

\bibitem[Bollob\'as et~al.(1998)Bollob\'as, Borgs, Chayes, and
  Kim]{bollobas-borgs-chayes-kim}
B\'ela Bollob\'as, Christian Borgs, Jennifer~T. Chayes, and Jeong~Han Kim,
  1998.
\newblock Lecture at the Workshop on the Interface between Statistical Physics
  and Computer Science, Torino, Italy.

\bibitem[Bollob\'as et~al.(2001)Bollob\'as, Borgs, Chayes, Kim, and
  Wilson]{bollobas-borgs-chayes-kim-wilson}
B\'ela Bollob\'as, Christian Borgs, Jennifer~T. Chayes, Jeong~Han Kim, and
  David~B. Wilson.
\newblock The scaling window of the 2-{SAT} transition.
\newblock {\em Random Structures \& Algorithms}, 18\penalty0 (3):\penalty0
  201--256, 2001.
\newblock arXiv:math.CO/9909031.

\bibitem[Chayes et~al.(1986)Chayes, Chayes, Fisher, and
  Spencer]{chayes-chayes-fisher-spencer:prl}
J.~T. Chayes, L.~Chayes, Daniel~S. Fisher, and T.~Spencer.
\newblock Finite-size scaling and correlation lengths for disordered systems.
\newblock {\em Physical Review Letters}, 57\penalty0 (24):\penalty0 2999--3002,
  1986.

\bibitem[Chayes et~al.(1989)Chayes, Chayes, Fisher, and
  Spencer]{chayes-chayes-fisher-spencer:cmp}
J.~T. Chayes, L.~Chayes, Daniel~S. Fisher, and T.~Spencer.
\newblock Correlation length bounds for disordered {Ising} ferromagnets.
\newblock {\em Communications in Mathematical Physics}, 120\penalty0
  (3):\penalty0 501--523, 1989.

\bibitem[Chv\'atal and Reed(1992)]{chvatal-reed:sat}
V.~Chv\'atal and B.~Reed.
\newblock Mick gets some (the odds are on his side).
\newblock In {\em Proceedings of the 33rd Symposium on the Foundations of
  Computer Science}, pages 620--627, 1992.

\bibitem[Crawford and Auton(1993)]{crawford-auton:sat}
James~M. Crawford and Larry~D. Auton.
\newblock Experimental results on the crossover point in satisfiability
  problems.
\newblock In {\em Eleventh National Conference on Artificial Intelligence
  (AAAI-93)}, pages 21--27, 1993.

\bibitem[Crawford and Auton(1996)]{crawford-auton:3-sat}
James~M. Crawford and Larry~D. Auton.
\newblock Experimental results on the crossover point in random 3{SAT}.
\newblock {\em Artificial Intelligence}, 81\penalty0 (1--2):\penalty0 59--80,
  1996.

\bibitem[Dubois et~al.(2000)Dubois, Boufkhad, and
  Mandler]{dubois-boufkhad-mandler:3-sat}
Olivier Dubois, Yacine Boufkhad, and Jacques Mandler.
\newblock Typical random 3-{SAT} formulae and the satisfiability threshold.
\newblock In {\em Proceedings of the Eleventh Annual ACM-SIAM Symposium on
  Discrete Algorithms}, pages 126--127, 2000.

\bibitem[Friedgut(1999)]{friedgut:sharp}
Ehud Friedgut.
\newblock Sharp thresholds of graph properties, and the $k$-{SAT} problem.
\newblock {\em Journal of the American Mathematical Society}, 12\penalty0
  (4):\penalty0 1017--1054, 1999.
\newblock With an appendix by Jean Bourgain.

\bibitem[Gent and Walsh(1994)]{gent-walsh:sat}
Ian~P. Gent and Toby Walsh.
\newblock The {SAT} phase transition.
\newblock In A.~Cohn, editor, {\em 11th European Conference on Artificial
  Intelligence}, pages 105--109. John Wiley \& Sons, Ltd., 1994.

\bibitem[Grimmett(1989)]{grimmett:percolation}
Geoffrey Grimmett.
\newblock {\em Percolation}.
\newblock Springer-Verlag, 1989.

\bibitem[Hayes(1997)]{hayes:sat}
Brian Hayes.
\newblock Can't get no satisfaction.
\newblock {\em American Scientist}, 85\penalty0 (2):\penalty0 108--112, 1997.

\bibitem[Kaporis et~al.(2002)Kaporis, Kirousis, and
  Lalas]{kaporis-kirousis-lalas:3.42}
Alexis~C. Kaporis, Lefteris~M. Kirousis, and Efthimios~G. Lalas.
\newblock The probabilistic analysis of a greedy satisfiability algorithm,
  2002.
\newblock Manuscript, presented at the European Symposium of Algorithms.

\bibitem[Kaporis et~al.(2000)Kaporis, Kirousis, and
  Stamatiou]{kaporis-kirousis-stamatiou:coloring}
Alexis~C. Kaporis, Lefteris~M. Kirousis, and Yannis~C. Stamatiou.
\newblock A note on the non-colorability threshold of a random graph.
\newblock {\em Electronic Journal of Combinatorics}, 7\penalty0 (1), 2000.
\newblock Paper \#R29.

\bibitem[Kaporis et~al.(2001)Kaporis, Kirousis, Stamatiou, Vamvakari, and
  Zito]{kaporis-kirousis-stamatiou-vamvakari-zito:sat}
Alexis~C. Kaporis, Lefteris~M. Kirousis, Yannis~C. Stamatiou, Malvina
  Vamvakari, and Michele Zito.
\newblock The unsatisfiability threshold revisited, 2001.
\newblock Manuscript.

\bibitem[Kirkpatrick(2000)]{kirkpatrick:personal}
Scott Kirkpatrick, 2000.
\newblock Personal communication.

\bibitem[Kirkpatrick et~al.(1993)Kirkpatrick, {Gy\"orgyi}, Tishby, and
  Troyansky]{kirkpatrick-gyorgyi-tishby-troyansky:k-sat}
Scott Kirkpatrick, G\'eza {Gy\"orgyi}, Naftali Tishby, and Lidror Troyansky.
\newblock The statistical mechanics of $k$-satisfaction.
\newblock In Jack~D. Cowan, Gerald Tesauro, and Joshua Alspector, editors, {\em
  Advances in Neural Information Processing Systems}, volume~6, pages 439--446.
  Morgan Kaufmann Publishers, 1993.

\bibitem[Kirkpatrick and Selman(1994)]{kirkpatrick-selman:critical}
Scott Kirkpatrick and Bart Selman.
\newblock Critical behavior in the satisfiability of random {Boolean}
  expressions.
\newblock {\em Science}, 264:\penalty0 1297--1301, 1994.

\bibitem[{\L uczak}(1991)]{luczak:chromatic}
Tomasz {\L uczak}.
\newblock The chromatic number of random graphs.
\newblock {\em Combinatorica}, 11\penalty0 (1):\penalty0 45--54, 1991.

\bibitem[Mitchell et~al.(1992)Mitchell, Selman, and
  Levesque]{mitchell-selman-levesque}
D.~Mitchell, B.~Selman, and H.~Levesque.
\newblock Hard and easy distributions of {SAT} problems.
\newblock In {\em Proc.\ 10th National Conference on Artificial Intelligence},
  pages 459--465, 1992.

\bibitem[Monasson and Zecchina(1998)]{monasson-zecchina:replica}
R\'emi Monasson and Riccardo Zecchina, 1998.
\newblock Informal talk on the replica method and 2-{SAT}.

\bibitem[Monasson et~al.(1999{\natexlab{a}})Monasson, Zecchina, Kirkpatrick,
  Selman, and Troyansky]{monasson-zecchina-kirkpatrick-selman-troyansky:rsa}
R\'emi Monasson, Riccardo Zecchina, Scott Kirkpatrick, Bart Selman, and Lidror
  Troyansky.
\newblock $2+p$-{SAT}: Relation of typical-case complexity to the nature of the
  phase transition.
\newblock {\em Random Structures \& Algorithms}, 15\penalty0 (3 and
  4):\penalty0 414--435, 1999{\natexlab{a}}.

\bibitem[Monasson et~al.(1999{\natexlab{b}})Monasson, Zecchina, Kirkpatrick,
  Selman, and Troyansky]{monasson-zecchina-kirkpatrick-selman-troyansky:nature}
R\'emi Monasson, Riccardo Zecchina, Scott Kirkpatrick, Bart Selman, and Lidror
  Troyansky.
\newblock Determining computational complexity from characteristic `phase
  transitions'.
\newblock {\em Nature}, 400:\penalty0 133--137, 1999{\natexlab{b}}.

\bibitem[Pittel et~al.(1996)Pittel, Spencer, and
  Wormald]{pittel-spencer-wormald:k-core}
Boris Pittel, Joel Spencer, and Nicholas Wormald.
\newblock Sudden emergence of a giant $k$-core in a random graph.
\newblock {\em Journal of Combinatorial Theory, Series B}, 67\penalty0
  (1):\penalty0 111--151, 1996.

\bibitem[Selman and Kirkpatrick(1996)]{selman-kirkpatrick:cost}
Bart Selman and Scott Kirkpatrick.
\newblock Critical behavior in the computational cost of satisfiability
  testing.
\newblock {\em Artificial Intelligence}, 81\penalty0 (1--2):\penalty0 273--295,
  1996.

\end{thebibliography}
\bibliographystyle{plainnat}

\end{document}